\documentclass[a4paper,12pt]{article}

\usepackage{amsfonts}
\usepackage{amscd,color}
\usepackage{amsmath,amsfonts,amssymb,amscd}
\usepackage{indentfirst,graphicx,epsfig}
\usepackage{graphicx}
\input{epsf}
\usepackage{graphicx}
\usepackage{epstopdf}
\usepackage{caption}

\newcommand{\ml}{l\kern-0.55mm\char39\kern-0.3mm}

\title{\textbf{Graph colorings under global structural conditions\footnote{Supported by
NSFC No.11871034 and 11531011. $^{**}$Corresponding author.}}}
\author{\small Xuqing Bai, Xueliang Li$^{**}$\\
\small Center for Combinatorics and LPMC\\
\small Nankai University, Tianjin 300071, China\\
\small Email: baixuqing0@163.com, lxl@nankai.edu.cn\\
}
\date{}
\begin{document}
\maketitle
\begin{abstract}
More than ten years ago in 2008, a new kind of graph coloring appeared in graph theory,
which is the {\it rainbow connection coloring} of graphs, and then followed
by some other new concepts of graph colorings, such as {\it proper connection coloring,
monochromatic connection coloring, and conflict-free connection coloring} of graphs, etc.
In about ten years of our consistent study, we found that these new concepts of graph
colorings are actually quite different from the classic graph colorings.
These {\it colored connection colorings} of graphs are brand-new colorings
and they need to take care of global structural properties
(for example, connectivity) of a graph under the colorings;
while the traditional colorings of graphs are colorings under which only local structural
properties (adjacent vertices or edges) of a graph are taken care of.
Both classic colorings and the new colored connection colorings can produce the so-called chromatic numbers.
We call the colored connection numbers the {\it global chromatic numbers}, and the classic or traditional
chromatic numbers the {\it local chromatic numbers}. This paper intends to clarify the difference
between the colored connection colorings and the traditional colorings, and finally to propose
the new concepts of global colorings under which global structural properties of the colored graph are kept,
and the global chromatic numbers.\\[3mm]
{\bf AMS Classification 2020:} 05C15, 05C40.

\end{abstract}

\section{Introduction}

The classic or traditional coloring theory is a very important and active
field in graph theory, and has invented many very powerful tools, methods, techniques,
and very beautiful results and theory. The most famous one that should be mentioned is
the so-called Four-Color Theorem (or Problem). More than ten years ago, a new kind of
graph coloring appeared, namely, the rainbow connection coloring of graphs,
which was introduced by Chartrand et al. \cite{CJMZ1} in 2008, and then followed by
some other new concepts of colorings, such as proper connection coloring \cite{ALLZ, BFGMMMT},
monochromatic connection coloring \cite{CY}, and conflict-free connection coloring \cite{CJV}, etc. \cite{JBY}.
For more information about them, we refer to \cite{CHLM, LM, LMQ, LSS, LS1, LS2, LW}.
We have been devoting to the study of these new colorings in recent years,
but there is always a question that was asked by ourselves or others, that is, does the
study of these new colorings make significant sense ? In other words,
are these colorings merely introduced as some new concepts in graph coloring theory ?
Or, are they really substantially different objects ? What is the significant difference
between these new colorings and the classic colorings ? As a result of meditation,
we finally realized that they really have quite different flavors. Colored connection
colorings are colorings under which global structural properties of graphs are kept.
However, traditional colorings are colorings under which only local structural properties
are guaranteed. Therefore, these colored connection colorings are called {\it global colorings}
of graphs and the corresponding chromatic numbers are called {\it global chromatic numbers};
while the classic or traditional colorings are called {\it local colorings} of graphs
and the corresponding chromatic numbers are called {\it local chromatic numbers}.

This paper aims to clarify the difference between the colored connection coloring
theory and the classic or traditional coloring theory. As a result, the new concepts of
global colorings and global chromatic numbers of graphs are proposed.

For any notation or terminology not defined here, we refer to Bondy and Murty \cite{BM}.

\section{Global structural properties vs Local structural properties}

Connectivity and colorings of graphs are the two most popular and important
concepts in graph theory. As a beginner of a graph theorist, one
met first are these two basic concepts. Let us recall them first.

To say that a graph is connected means that there is a path
between any two vertices in the graph to connect them,
which means that the graph is one ``block" as a whole. We know that the
local connectivity of a graph cannot guarantee that the entire
graph is connected, and vice versa, obviously.
In other words, connectedness is not only concerned with the local structure
of a vertex, but also the connection relationship between a vertex and
any other vertex in the graph. So, connectivity cares about global structures
of a graph.

We now turn to the colorings of graphs,
which are generally classified as edge colorings and vertex colorings,
and proper colorings are usually studied.
For edge colorings, a {\it proper edge coloring} of a graph is a
coloring to the edges of the graph such that the colors of the edges incident with
one vertex are different from each other, and the colors of the other edges
are not required. In other words, the edges it cares
about are only in the local structure of a vertex.
So, traditional colorings only care about local structures of a graph.

\section{ When connectivity meets colorings}

Given a colored graph, without loss of generality, we consider an
edge colored graph, how do we require connectivity in the colored graph ?
Obviously, we cannot ignore the colors that are used on the graph.
So, what should the concept of connectivity that takes the coloring
into account look like ? This is the problem of colored
connectivity that we are studying now. We know that connectivity is
defined by paths, and so we have to use colored paths now,
that is, a path with a certain color pattern or property.
If this pattern or property is denoted by $\mathcal{P}$,
then this kind of paths are also called {\it $\mathcal{P}$-paths}.

An edge colored graph $G$ is called {\it $\mathcal{P}$-connected},
if between any two vertices of $G$  there is a $\mathcal{P}$-path
to connect them. This edge coloring is also called a {\it $\mathcal{P}$-connection coloring}
of the graph. Given an uncolored connected graph $G$, the
{\it $\mathcal{P}$-connection coloring number}, denoted by $\mathcal{P}c(G)$,
is the {\it minimum (maximum)} number of colors that are needed in order to
make $G$ $\mathcal{P}$-connected. When the property $\mathcal{P}$ is {\it rainbow or proper},
we take the {\it minimum number} of colors and call it {\it rainbow connection number} or
{\it proper connected number}, denoted by $rc(G)$ and $pc(G)$, respectively;
When the property $\mathcal{P}$ is {\it monochromatic},
we take the {\it maximum number} of colors (because taking minimum number of colors is meaningless),
and call it the {\it monochromatic connected number}, denoted by $mc(G)$.
Of course, corresponding to the previous terms, the graph is then called
{\it rainbow (proper, monochromatic) connected}, the coloring is called
{\it rainbow (proper, monochromatic) connection coloring}, and
the chromatic number is called {\it rainbow (proper, monochromatic) connection number}.

We also call a $\mathcal{P}$-connection coloring a {\it colored connection coloring}.
How to choose the property or pattern $\mathcal{P}$ ? In practical applications,
it depends on the needs of customers, such as the security settings of communications, interference-free
transmission of signals, or selection of aviation lines, etc. There have been quite a lot of applications of
examples in existing literature. In theoretical research, it depends on researchers' own interests,
for which there have also appeared many new concepts of colored connection colorings of graphs in separate literature.

One thing we would like to say is that colored connection colorings belong to the topic of graph colorings. This is
because they mainly concern how to color a graph by using minimum (maximum) numbers of colors, i.e., the optimal
colorings and the chromatic numbers. They have close relations with graph connectivity, but not substantial relations,
since graphs with the same connectivity may have quite different colored connection colorings (numbers).

\section{Compared with traditional colorings: Local chromatic numbers vs Global chromatic numbers}

From the above definition, it can be seen that each colored
connection number is a kind of ``chromatic number". It is the
same as the classic chromatic number, that is, given a graph
and we color it under the restriction of some structural properties
of the graph and find the minimum (or maximum) number of colors needed.
However, the colored connection number needs to consider the
connectivity of the graph, and therefore, as we stated before, the
colored connection number is a chromatic number under the condition
of global structures, referred to as a {\it global chromatic number};
while a traditional or classic chromatic number only needs to
satisfy some local structural conditions (for edge coloring, the colors
only on adjacent edges are completely distinguishable), that is,
a chromatic number under local structural conditions,
referred to as a {\it local chromatic number}, for short.

As we know, classic coloring theory has developed many powerful
methods, tools, techniques, and beautiful results and theory.
For example, we can use combinatorial enumeration method -- introducing chromatic polynomials \cite{B};
probabilistic method -- using the famous Local Lemma \cite{EL}; or algebraic
method -- Combinatorial Nullstellensatz \cite{A}, and so on.
Can we use these methods that are very powerful in classic coloring theory to
the study of colored connection colorings ? Let us examine them one by one in the following subsections.

\subsection{Chromatic polynomial - A combinatorial enumeration method}

Can we have a colored connection coloring polynomial ?

Let us recall the definition of the traditional chromatic polynomial.
Given a positive integer $k$ and a graph $G$,
count the number of how many ways to color the graph by proper $k$-colorings,
and get an enumeration function $f(G,k)$. It uses two recursive operations:
deletion and contraction of edges to prove that the function $f(G,k)$ is equal to
$f(G-e, k) - f(G \cdot e, k)$ \cite{BM}, and eventually it is a polynomial in $k$,
and so it is called the chromatic polynomial of $G$. Of course this is obtained for
vertex colorings of a graph. One can also get an analogous polynomial for edge colorings
of a graph, simply by just considering the line graph $L(G)$ of $G$.
Thus, the Four-Color Conjecture of a planar graph
is that for any planar graph $G$, we always have $f(G,4)>0$,
that is, there always is a way to properly 4-color any planar graph.
Why is there such a good polynomial ? This is the benefit that classic
colorings only involve guaranteeing some local structural properties of
a vertex.

Now, going back to the problem we started with in this section,
we want a polynomial for colored connection colorings, or to be exact,
we want an enumeration function for counting the number of ways of
colored connection colorings to a given graph.
So, to answer the question, we first need to clarify what kind of coloring
is a $\mathcal{P}$-connection coloring for a given property $\mathcal{P}$ ?
Let us analyze this problem as follows. First, we have the choice of
$n\choose 2$ vertex pairs; and then for a given pair of vertices,
there is a $\mathcal{P}$-path between them,
but there are a lot of paths between them (could be exponential many in general).
Besides, these paths can be distributed anywhere in the entire graph,
which makes the problem extremely complicated.
How can we enumerate our so-called ``proper" colorings ? For general graphs,
it is just too many to count. So, this should be the difficulty that
the $\mathcal{P}$-connection colorings have.

However, to determine whether an edge coloring of a graph is proper
should be easy thing to do, just check each vertex to see if the edges
incident with it are colored differently. This is a polynomial time problem.
Nevertheless, given an edge coloring of a graph, check whether
it is a rainbow connection coloring is NP-complete \cite{CFMY}.
Therefore, we can hardly expect to have a rainbow connection coloring
polynomial, or even if an exact expression of the enumeration function of
rainbow connection colorings.

\subsection{Local Lemma - A probability method}

We know that the Local Lemma is a powerful tool in
probabilistic methods. Let us recall the following form of the Local Lemma as a symmetric form.

\noindent {\bf Theorem:} {\upshape\cite{EL}}
Let $A_i$, $i\in N$, be events in a probability space $(\Omega,P)$
having a dependency graph with maximum degree $d$.
Suppose that $P(A_i)\leq 1/(e(d + 1)),$ $i\in N$.
Then  $P(\cap_{i\in N} \overline{A_i})>0$.

For traditional colorings, say proper edge coloring.
If the edges incident with each vertex in a graph $G$ are colored in
different colors, then the coloring is a proper edge coloring
of $G$. From the conditions of the Local Lemma,
we first need to construct a bad event. In a coloring of $G$,
we use $A_v$ to represent a bad event, that is, there are two edges incident
with a vertex $v$ in the graph $G$ that have the same color. So,
$\overline{A_v}$ is a good event. Hence, if $P(\cap_{v\in V(G)} \overline{A_v})>0$,
then we know that there exists a proper edge coloring for $G$.
In traditional coloring, the bad event $A_v$ is well described,
because we only need to consider the colors of the edges incident with
the vertex $v$. This is the benefit of the classic colorings,
involving only local structural (adjacent edges) properties of a vertex.

However, for colored connection colorings, our bad event $A_{u,v}$
for a pair of vertices $u$ and $v$ means that there is no $\mathcal{P}$-path to
connect $u$ and $v$ in $G$. As one can image, there are many paths between $u$
and $v$ in the graph $G$, and they are intertwined, their lengths vary, and
they may involve all vertices and edges in the entire graph.
This makes it impossible for us to describe the bad event $A_{u,v}$ clearly and accurately.
As a result, we can hardly use Local Lemma, as such a good tool.

\subsection{Combinatorial Nullstellenensatz - An algebraic method}

The Combinatorial Nullstellenensatz is a generalization of one
variable polynomials with the most order roots.
Using this tool means that we can construct a polynomial expression
in combination with the graph property $\mathcal{P}$,
so that if the polynomial exists a non-zero solution,
then the graph has a coloring that satisfies the property $\mathcal{P}$.
Otherwise, the graph does not have a coloring to satisfy the property
$\mathcal{P}$.

In proper edge coloring, we can set each edge $e$ of the
graph as a variable $x_e$. For each vertex $v$ with degree $d(v)\geq 2$,
set the product of the differences of the variables corresponding to
pairs of the edges incident with $v$ to be a polynomial.
Then, multiply the polynomials corresponding to all vertices with degree
$d(v)\geq 2$ to obtain a multi-variate polynomial about proper coloring.
Now, assign values to the edge variables. If there is a set of assignments
to make the value of the polynomial non-zero, it means that the graph has a
proper edge coloring. Otherwise, if one of the items is zero, the coloring corresponding to
these assignments is not a proper edge coloring. This is a very clear description,
because in proper edge coloring, we only need to consider the colors
of the edges incident with some vertex, and the edges incident with a vertex
are clearly given, and the number of them is very limited.
Under this local structural restriction, we can easily construct
a polynomial expression to decide a proper edge coloring.

However, for example, for a rainbow connection coloring, that is, an
edge coloring of a graph $G$ so that there is a rainbow path in $G$ connecting
each pair of vertices of $G$, we have to consider all the paths between a given pair
of vertices. For a general graph, there are a large number of these paths,
which can be exponentially large, and the interaction of these paths is complicated,
and they even cover all edges of $G$. So, it is very complicated to express
all these possible paths clearly and accurately,
which makes us unable to succinctly give a polynomial expression to decide whether
the coloring is rainbow connected. This makes us to hardly use the beautiful
Combinatorial Nullstellenensatz.

{\it Discharge Technique} \cite{BM} is also popularly used in traditional vertex(edge) colorings. However,
it is not difficult to analyze that it can hardly be used in colored connection colorings.

It is worth mentioning that there is a paper \cite{AD} that tries to use the
Combinatorial Nullstellenensatz method to study the rainbow connectivity of graphs.
However, as one can see there, it is quite complicated and not practical.
Because only when the rainbow connection number $k$ is a constant, we can convert
the graph instance into a polynomial equation system in the time polynomial of $n$.

As a generalization of colored connection colorings, we then define the
following concept of colored $k$-connection colorings.

\section{ Colored $k$-edge(vertex) connection colorings}

In this section, we will introduce the concept of colored connectivity.

Let $k$ be a positive integer, an edge colored graph is called
$k$-edge(vertex) connected, if there are $k$ paths with disjoint edges
(internal vertices) between any two vertices in $G$ to connect them.
This edge(vertex) coloring is called a $k$-edge(vertex) connection
coloring of $G$. Given an uncolored graph, its $k$-edge(vertex)
$\mathcal{P}$-connection number, denoted by $\mathcal{P}c_k(G)$ is defined as
the minimum (maximum) number of colors required in a $k$-edge(vertex)
$\mathcal{P}$-connection coloring of $G$. Similarly, if $\mathcal{P}$ is
rainbow or proper, then we take the minimum number of colors;
while if $\mathcal{P}$ is monochromatic, we take the maximum number of colors.
Obviously, when $k=1$, the $k$-edge(vertex) $\mathcal{P}$-connection number and
$\mathcal{P}$-connection number are the same, which is called the $\mathcal{P}$-connection number of $G$,
simply denoted by $\mathcal{P}c(G)$. For examples, we have $rc(G)$, $pc(G)$, $cfc(G)$, and $mc(G)$
for the rainbow, proper, conflict-free, and monochromatic connection number, respectively.

At present, most of the literature are on the case $k=1$ (see
\cite{CLW,LSS,LS1,LS2} for examples and the references therein), and there are some results
on the case $k=2$ (see \cite{CJMZ2,CLL,FLM,LL2,LS3} for examples),
but there are very few or almost no results on the case $k=3$.
It is not that one does not like to study it, but because it is simply too difficult
to get reasonable or good results.

\section{Colored disconnection colorings}

Connectedness and disconnectedness of graphs are the two sides of one coin on
a graph. One can use paths to describe the connectivity of a graph, and one also
can use edge(vertex) cuts to do the same thing.
Given a graph $G$, for any two vertices $x$ and $y$ of $G$,
the maximum number of pairwise edge-disjoint $xy$-paths is
equal to the minimum number of edges in an $xy$-cut.
This is the statement of the famous Menger's Theorem.
Thus, we can also consider the edge connectivity of $G$ by using edge cuts.
Let $\mathcal{P}$ be a property or a pattern. An edge cut $R$ of an edge colored
connected graph $G$ is called a {\it $\mathcal{P}$-cut} if the edges in $R$ have
the property $\mathcal{P}$. Let $u$ and $v$ be two vertices of $G$.
A $u$-$v$-\emph{$\mathcal{P}$-cut} is a $\mathcal{P}$-cut $R$ of $G$
such that $u$ and $v$ belong to different components of $G-R$.
An edge colored graph $G$ is called \emph{$\mathcal{P}$-disconnected}
if for every two vertices $u$ and $v$ of $G$, there exists a
$u$-$v$-$\mathcal{P}$-cut in $G$ separating them.
In this case, the edge coloring is called a \emph{$\mathcal{P}$-disconnection coloring} of $G$.
For an uncolored connected graph $G$, the \emph{$\mathcal{P}$-disconnection number}
of $G$, denoted by $\mathcal{P}$d$(G)$, is the minimum (maximum) number
of colors that are needed in order to make $G$ $\mathcal{P}$-disconnected.
When $\mathcal{P}$ is rainbow or proper, the $\mathcal{P}$-disconnection number takes the
minimum number of colors, denoted by $rd(G)$ or $pd(G)$;
when $\mathcal{P}$ is monochromatic, the $\mathcal{P}$-disconnection number
takes the maximum number of colors, denoted by $md(G)$.
By Menger's Theorem, colored disconnection also studies the
colored connection of a graph $G$, which is the two ways of one thing, the colored connection.
Obviously, the colored disconnection number $\mathcal{P}d(G)$ is also a global chromatic number,
and there have been some results about it; see \cite{BCL,CDHHZ,CLLW,PL1}.

Since we already have colored connection (by colored paths) and disconnection (by colored edge cuts) in a graph,
it is worth to raise a question: How to generalize the famous Menger's theorem to a colored version ?
At the moment, we have not found a correct statement of a colored version for the Menger's theorem. It would be
very interesting to create and prove a colored version of the Menger's theorem.

In addition to the connectivity, we know that graphs have many other
global structural properties, such as spanning trees, Hamiltonian cycles (paths) and so on.
Similarly, combining them with graph colorings, one may think about introducing some other new
global colorings and global chromatic numbers, for example, proper-walk connection coloring and number
in a recent paper \cite{JBY}.

\section{Other versions of global colorings}

The previous discussions are all about the study of edge colored
version of graphs, and it is natural to consider the study of other
forms of colored versions. So far, there have been results on
vertex colored version and total colored (both vertex and edge colored) version,
as well as strong version (using shortest paths).
There are also studies on colored version of digraphs, and even colored version of hypergraphs.
Furthermore, one can also consider to carry out study on colored connection and disconnection
for matroids. The definitions of all of them can be elaborated out naturally,
and so their details are omitted here.

It is worth mentioning that there are some papers \cite{BJJ, JLLM} studying proper rainbow connection
colorings and disconnection colorings of graphs, which means that the graph is properly
colored and under this very coloring the graph is also rainbow connected or disconnected.
This kind of colorings can be viewed as colorings under both local and global structural conditions,
and as a whole, they are still global colorings of graphs. One may think about list colorings combined with
rainbow connection colorings and adjacency edge(vertex) distinguishable colorings \cite{CWZ} combined with rainbow
connection colorings, to get {\it list rainbow connection colorings} and {\it adjacency edge(vertex) distinguishable
rainbow connection colorings}, etc.

It is also worth mentioning that there is a research topic closely related
to this, which is about the existence of colored subgraphs (substructures)
in colored graphs. Given an edge colored graph $G$, under what kind of color conditions,
there are rainbow (proper, monochromatic, etc.)
substructures with $k$ edges, such as cycles, paths, trees, matchings, and so on.
When $k$ is large enough, the corresponding substructures are rainbow (proper, monochromatic, etc.)
Hamiltonian cycles, Hamiltonian paths, spanning trees, perfect matchings, etc.
In general, there are color degree conditions, forbidden subgraph conditions, and so on.
There has been a large amount of published literature in this subject; see \cite{AFR,CKRY,CHH,CKW,FNXZ,KL,LNXZ,LNZ,L1,L2}
for examples and the references therein. But, this kind of research is obviously different from the one here.
It is the {\it theory of colored graphs}, while the study here belongs to the {\it theory of graph colorings}.
The former is to study the existence of colored substructures in colored graphs;
while the latter is to study the graph colorings of an uncolored graph under certain structural properties,
requiring the minimum or maximum number of colors, that is, the chromatic number.

\section{Concluding remarks: Global colorings and global chromatic numbers}

In the above sections, we first recalled a global structural property of graphs, i.e.,
the connectivity, and a graph object and parameter defined by local structural conditions of graphs,
i.e., the edge(vertex) coloring and chromatic index(number), which we called it a local coloring and a local chromatic number.
Then, when connectivity meets graph colorings, we introduced the concept of colored connection colorings and numbers, which are
graph objects and parameters defined by global structural conditions of graphs, we called them global colorings and global chromatic numbers.
To understand the significant difference between the global and local
colorings and chromatic numbers, we employed three examples of popularly used methods to analyze
the reasons why the local one can use them successfully but the global ones can hardly do.

For further study along with this topic, we will introduce a general setting of global colorings
and global chromatic numbers as follows.

\noindent {\bf Definition:} Let $G$ be an edge colored graph, and $\mathcal{P}$ be a graph property.
If $G$ has the property $\mathcal{P}$, we say that $G$ is a {\it $\mathcal{P}$-graph}, and the edge coloring is
called a {\it $\mathcal{P}$-coloring} of $G$.
For examples, $\mathcal{P}$ can be the global property: {\it rainbow (proper, monochromatic, etc.)
$k$-connected ($k\geq 1$)}, or {\it rainbow (proper, monochromatic, etc.)
disconnected}, or {\it proper rainbow connected or disconnected}, etc. Then, we have that $G$ is a {\it rainbow (proper, monochromatic, etc.)
$k$-connected graph}, or a {\it rainbow (proper, monochromatic, etc.)
disconnected graph}, or a {\it proper rainbow connected or disconnected graph}, etc. The corresponding edge coloring is then called a {\it rainbow
(proper, monochromatic, etc.) $k$-connection coloring}, or a {\it rainbow (proper, monochromatic, etc.)
disconnection coloring}, or a {\it proper rainbow connection or disconnection coloring}, etc. Given an uncolored graph $G$, we consider the minimum (maximum)
number of colors required to make $G$ have the property $\mathcal{P}$. This number is called the {$\mathcal{P}$-chromatic number
of $G$, denoted by $\chi_\mathcal{P}(G)$.
Similarly, when the property $\mathcal{P}$ is rainbow, or proper, or proper rainbow,
the minimum number of colors is taken; while when $\mathcal{P}$ is
monochromatic, the maximum number of colors is taken.

Based on this definition, it is expected that more and more interesting
new global colorings and chromatic numbers will appear. This, on the one hand,
depends on the needs of practical applications to stimulate introducing new concepts
and new theories; on the other hand, it depends on the interest of theoretical
researchers to introduce new concepts with theoretical significance.
Compared with the traditional coloring theory, this is a completely
different field of study. Global colorings belong to a brand-new coloring theory of graphs.
To clarify this viewpoint, as we exemplified in Section 4 above,
the first thing to do is to grasp what coloring is a traditional
proper edge(vertex) coloring, and what coloring is a rainbow (proper, monochromatic)
connection coloring, as we discussed before ?  What is the complexity
brought about by their respective structural condition requirements, and what kind of conditions should
be checked ? Is it clear and succinct to check these conditions ?
What is the time complexity to check each condition ?
For this kind of new colorings, we need to invent new methods and tools to deal with them,
and eventually to establish a new system of theory.
More and more efforts have to be paid to make it perfect.

\end{document}